\documentclass[10pt,twoside]{article}
\usepackage{graphicx}
\usepackage{amsmath}
\usepackage{Latex-document,amsfonts,amsmath,amssymb}
\newtheorem{thm}{Theorem}[section]

\newtheorem{lem}[thm]{Lemma}

\newtheorem{defn}[thm]{Definition}
\newtheorem{prop}[thm]{Proposition}

\newtheorem{axiom}{Axiom}
\newtheorem{3+}{Axiom $3^+$}
\newtheorem{example}{Example}
\newtheorem{example*}{Example}

\newtheorem{5'}{Axiom B$5'$}

\newcommand{\wdist}[2]{\ensuremath{\lVert #1, #2 \rVert^\alpha_
{\mathrm{w}}}}

\newcommand{\sdist}[2]{\ensuremath{\lVert #1,#2 \rVert^\alpha_{\mathrm{s}}}}
\newcommand{\distd}[2]{\ensuremath{\lVert #1, #2 \rVert_{*}}}
\newcommand{\onedist}[2]{\ensuremath{\lVert #1,#2 \rVert_1}}
\renewcommand{\thesection}{\arabic{section}}
\renewcommand{\thesubsection}{\thesection.\arabic{subsection}.}
\def\addsec{\addtocounter{section}{1} \setcounter{thm}{0}}

\title{\bf Applications of Orbit Equivalence \vskip -2mm
to Actions of Discrete Amenable Groups\vskip 6mm}

\author{Daniel J. Rudolph\vspace*{-0.5cm}\thanks{Department of Mathematics,
University of Maryland, College Park, MD 20742, USA. E-mail:
djr@math.umd.edu }}
\date{\vspace{-8mm}}

\begin{document}

\maketitle

\thispagestyle{first} \setcounter{page}{339}

\begin{abstract}

\vskip 3mm

Since the work of Ornstein and Weiss in 1987 ({\bf Entropy and isomorphism theorems for actions of amenable
groups}, {\it J. Analyse Math.}, {\bf 48} (1987)) it has been understood that the natural category for classical
ergodic theory would be probability measure preserving actions of discrete amenable groups.  A conclusion of this
work is that all such actions on nonatomic Lebesgue probability spaces were orbit equivalent.  From this
foundation two broad developements have been built. First, a full generalization of the various equivalence
theories, including Ornstein's isomorphism theorem itself, exists.  Fixing the amenable group $G$ and an action of
it, one can define a metric-like notion on the full-group of the action, called a size. A size breaks the orbit
equivalence class of a single action into subsets, those reachable by a Cauchy sequence (in the size) of full
group perturbations.  These subsets are the equivalence classes associated with the size.  Each size possesses a
distinguised ``most random" set of classes, the ``Bernoulli" classes of the relation.  An Ornstein-type theorem
can be obtained.  Many naturally occuring equivalence relations can be described in this way.  Perhaps most
interesting, entropy itself can be so described.
Second, one can use the characterization of discrete amenable
actions as those which are orbit equivalent to a action of ${\mathbb Z}$ to lift theorems from actions of $\mathbb
Z$ to those of arbitrary amenable groups.  The most interesting of these are first, that actions of completely
positive entropy (called $K$-systems for ${\mathbb Z}$ actions) are mixing of all orders (proven jointly with B.
Weiss) and that such actions have countable Haar spectrum (proven by Golodets and Dooley).  As all ergodic actions
are orbit equivalent, only ergodicity is preserved by orbit equivalences in general, but by considering orbit
equivalences restricted to be measurable with respect to a sub-$\sigma$ algebra, many properties relative to that
algebra are preserved.  This provides the tool for this method to succeed. \vskip 4.5mm

\noindent {\bf 2000 Mathematics Subject Classification:} 28D15,
37A35.

\noindent {\bf Keywords and Phrases:} Amenable group, Orbit
equivalence, Entropy.
\end{abstract}

\vskip 12mm

\section*{1. Definitions and examples of sizes} \addsec
%\label{section 1}\setzero
\vskip-5mm \hspace{5mm}

Our goal in this section is to describe a metric-like notion on
the full group of a measure preserving action of an amenable
group and show how this leads to various restricted orbit
equivalence theories. This work can be found in complete detail
in {\bf Restricted Orbit Equivalence for Actions of Discrete
Amenable Groups} by D.J. Rudolph and J. Kammeyer, {\it Cambridge
Tracts in Mathematics} \# 146.

Let $(X,\mathcal{F},\mu)$ be a fixed nonatomic Lebesgue
probability space.  Let $G$ be an infinite discrete amenable
group. Let $\mathcal{O}\subseteq X\times X$ be an ergodic, measure
preserving, hyperfinite equivalence relation.  For our purposes,
this simply means that $\mathcal{O} = \{(x,T_g(x))\}_{g\in G}$
where $T:G\times X\to G$ (written of course $T_g(x)$)  is some
ergodic and free, measure preserving action of $G$ on $X$.

\begin{defn}Let $G$ be an infinite countable discrete amenable group.  A
$\mathbf{G}${\bf -arrangement} $\alpha$ is any map from
$\mathcal{O}$ to $G$ that satisfies:

\begin{itemize}
\item[(i)] $\alpha$ is 1-1 and onto, in that for a.e. $x\in X$,
for all $g\in G$, there is a unique $x'\in X$ with
$\alpha(x,x')=g$.  We write $x'=T_g^{\alpha}(x)$;

\item[(ii)] $\alpha$ is measurable and measure
preserving, i.e. for all $A\in \mathcal{F}, g\in G$, both
$T_g^{\alpha}(A)\in \mathcal{F}$ and $\mu
(T_g^{\alpha}(A))=\mu(A)$; and

\item[(iii)] $\alpha$ satisfies the cocycle equation
$\alpha(x_2,x_3)\alpha(x_1,x_2)=\alpha(x_1,x_3)$.

\end{itemize}
As $G$ will not vary for our considerations we will abbreviate
this as an {\bf arrangement}.  Let $\mathcal{A}$ denote the set of
all such arrangements.
\end{defn}

\begin{lem}
$\alpha$ is a $G$-arrangement if and only if there is a measure
preserving ergodic free action of $G$, $T$, whose orbit relation
is $\mathcal{O}$\index{orbit relation} such that
$\alpha(x,T_g(x))=g$ for all $(x,T_g(x))\in\mathcal{O}$.
\end{lem}

Thus the vocabulary of $G$-arrangements on $\mathcal{O}$ is
precisely equivalent to the vocabulary of $G$-actions whose orbits
are $\mathcal{O}$.  For a $G$-arrangement $\alpha$, we write
$T^\alpha$ for the corresponding action.  For a $G$-action $T$, we
write $\alpha _T$ for the corresponding $G$-arrangement.

\begin{defn}
The {\bf full group} of $\mathcal{O}$ is the group (under
composition) $\Gamma$ of all measure preserving invertible maps
$\phi : X\to X$ such that for $\mu$-a.e. $x\in X$, $(x,\phi(x))\in
\mathcal{O}$.
\end{defn}

\begin{defn}
A $\mathbf{G}${\bf -rearrangement}of $\mathcal{O}$ is a pair $(\alpha,\phi)$, where $\alpha$ is a \linebreak
$G$-arrangement of $\mathcal{O}$ and $\phi \in \Gamma$.   As $G$ is fixed for our purposes we will abbreviate this
as a {\bf rearrangement}.  Let $\mathcal{Q}$ denote the set of all such rearrangements. \index{rearrangement(s) !
def}
\end{defn}

Intuitively, a rearrangement is simply a change (i.e.
rearrangement) of an orbit from the arrangement $\alpha$ to the
arrangement $\alpha\phi$, where $\alpha\phi(x,x') =
\alpha(\phi(x),\phi(x'))$.  One can formalize such a rearrangement
in three different ways.  Set $\mathcal{B}$ to be the set of
bijections of $G$ and $\mathcal{B}$ the subgroup of $\mathcal{G}$
fixing the identity.  Both are topologized via the product
topology on $G^G$.  Notice there is a homomorphism $\hat
H:\mathcal{B}\to\mathcal{G}$ given by $\hat
H(q)(g)=q(\mathrm{id})^{-1}q(g)$.   The kernel of $\hat H$ consist
of the left translation maps.

To a rearrangment we can associate a family of functions
$q^{\alpha,\phi}_x\in\mathcal{B}$ where
\[q^{\alpha,\phi}_x(g)=\alpha(x,\phi(T_g^\alpha(x))).\]\index{q@$q_x^{\alpha,\phi}$}
Now suppose $\alpha$ and $\beta$ are two arrangements of the
orbits $\mathcal{O}$.  Regard the first as an initial arrangement
and the second as a terminal arrangement.  We can associate to
this pair and any point $x$ a bijection from $G$ fixing the
identity that describes how the arrangement of the orbit has
changed:
\[h^{\alpha,\beta}_x(g)=\beta(x, T^{\alpha}_g(x)).\]
Notice here that $\hat
H(q_x^{\alpha,\phi})=h_x^{\alpha,\alpha\phi}.$

 Write
$h^{\alpha,\beta}:X\to \mathcal{G}.$

The third way to view a rearrangement pair has a symbolic dynamic
flavor.  For each orbit $\mathcal{O}(x)=\{x'; (x,x')\in
\mathcal{O}\}$, a rearrangement $(\alpha,\phi)$ also gives rise to
a natural map $G\to G$ (not a bijection though), given by
\index{f@$f_x^{\alpha,\phi}$}
\[f_x^{\alpha,\phi}(g)=\alpha(T^{\alpha}_g(x),\phi(T^{\alpha}_
g (x)).\]
Visually, regarding $\mathcal{O}(x)$ laid out by
$\alpha$ as a copy of $G$, $\phi$ translates the point at
position $g$ to position $f^{\alpha,\phi}_x(g)g$.

There is a natural link between the three functions
$h^{\alpha,\alpha\phi}$, $q^{\alpha,\phi}$ and $f^{\alpha,\phi}$
as follows.  For any map $f:G\to G$ we define
\[Q(f)(g)=f(g)g\text{ and}\]
\[H(f)(g)=f(g)gf(\mathrm{id})^{-1}.\] \index{h@$H(f)$}
It is an easy calculation that
\[H(f^{\alpha,\phi})=h^{\alpha,\alpha\phi}\text{ and }Q(f^{\alpha,\phi}=q^{\alpha,\phi}.\]
Let $\{F_i\}$ be a fixed F\o lner sequence for $G$.  We will
describe a number of concepts in terms of the $F_i$.

We now consider three pseudometrics on the set of rearrangements.
These all arise from natural topologies on functions $G\to G$. As
$G$ is countable the only reasonable topology is the discrete one,
using the discrete 0,1 valued metric.  This topologizes $G^G$ as a
metrizable space with the product topology. This is the weakest
topology for which the evaluations $g:f\to f(g)$ are continuous
functions. Notice that $H$ is a continuous map from $G^G$ to
itself and the map $h\to h^{-1}$ on $\mathcal{G}$ is continuous.

Define a metric $d$ on $\mathcal{G}$ as follows. List the elements
of $G$ as $\{g_1=\mathrm{id},g_2,\dots\}$ and let $d_0$ be the 0,1
valued metric on $G$.  Set
\[d(h_1,h_2)=\sum_i
[d_0(h_1(g_i),h_2(g_i))+d_0(h_1^{-1}(g_i),h_2^{-1}(g_i)))]
2^{-(i+1)}.\] Notice that if $h_1$, $h_2$, $h_1^{-1},$ and
$h_2^{-1}$ agree on $g_1,\dots,g_i$ then $d(h_1,h_2)\leq2^{-i}$.
On the other hand if $d(h_1,h_2)<2^{-i}$ then $h_1$, $h_2$ and
their inverses agree on this list of $i$ terms.

\begin{lem}
 The metric $d$ on $\mathcal{G}$ \index{g@$\mathcal{G}$} gives
the restricted product topology and makes $\mathcal{G}$ a complete
metric space.\end{lem}

We can use this to define a complete $L^1$ metric on arrangements:

\[\onedist{\alpha}{\beta}=\int
d(h^{\alpha,\beta},\mathrm{id})\,d\mu.\] As
$d(h_1,h_2)=d(h_2^{-1}h_1,\mathrm{id})$ and
$(h^{\alpha,\beta}_x)^{-1}=h^{\beta,\alpha}_x$ we see that this is
a metric.

We can also define a metric similar to $d$ on $G^G$ itself making
it a complete metric space by just taking half of the terms in
$d$:
\[d_1(f_1,f_2)=\sum_i d_0(f_1(g_i),f_2(g_i))2^{-i}.\]  This also
leads to an $L^1$ metric on $G^G$-valued functions on a measure
space:

\[\onedist{f_1}{f_2}=\int d_1(f_1,f_2)\,d\mu.\]

These two $L^1$ distances now give us two families of $L^1$
distances on the full-group, one a metric the other a
pseudometric, associated with an arrangement $\alpha$:

\[\wdist{\phi_1}{\phi_2}=\int d(h^{\alpha,\alpha\phi_1},h^{\alpha,\alpha\phi_2}) =\onedist{\alpha\phi_1}{\alpha\phi_2}
\]
and
\[\sdist{\phi_1}{\phi_2}=\int d_1(f^{\alpha,\phi_1},f^{\alpha,\phi_2})\,d\mu =
\onedist{f^{\alpha,\phi_1}}{f^{\alpha,\phi_2}}.\]

The {\bf weak} $L^1$ distance, $\wdist{\cdot}{\cdot}$, is only a
pseudometric but the {\bf strong} $L^1$ distance,
$\sdist{\cdot}{\cdot}$, is a metric.

To describe the weak*-distance  between two arrangements let
$G^\star=G\cup\{\star\}$ be the one point compactification of $G$.
Now $(G^\star)^G$ is a compact metric space and hence the Borel
probability measures on $(G^\star)^G$, which we write as
$\mathcal{M}_1(G^\star)$, are a compact and convex space in the
weak* topology. Let $D(\mu_1,\mu_2)$ be an explicit metric giving
this topology.

We define the distribution pseudometric
 between two
rearrangements by
\[\distd{(\alpha,\phi)}{(\beta,\psi)}=D((f^{\alpha,\phi})^*(\mu),
(f^{\beta,\psi})^*(\nu)).\]

We can combine the two $L^1$-metrics on arrangements and the full
group to define a product metric on rearrangements in the form

\[\lVert
(\alpha_1,\phi_1),(\alpha_2,\phi_2)\rVert_1=\onedist{\alpha_1}
{\alpha_2}+\mu(\{x:\phi_1(x)\neq \phi_2(x)\}).\]  We end this
Section by relating this complete $L^1$-metric on rearrangements
to the distribution pseudometric.

We now define the notion of a {\bf size} $m$ on rearrangements
$(\alpha,\phi)$ as a family of pseudometrics $m_\alpha$ on the
full-group satisfying some simple relations to the metrics and
pseudometrics we just defined.

A size is a function
\begin{displaymath}
m:\mathcal{Q} \to \mathbb{R}^+
\end{displaymath}
such that, if we write
\begin{displaymath}
m_{\alpha}(\phi_1,\phi_2) \underset{\mathrm{defn}}{=}
m(\alpha\phi_1,\phi_1^{-1}\phi_2),
\end{displaymath}
then $m$ satisfies the following three axioms.

\begin{axiom}
For each $\alpha\in\mathcal{A}$, $m_{\alpha}$ is a pseudometric on
$\Gamma$.
\end{axiom}

\begin{axiom}
For each $\alpha\in\mathcal{A}$, the identity map
\begin{displaymath}
(\Gamma,m_{\alpha}) \overset{\mathrm{id}}{\rightarrow} (\Gamma,
\wdist{\cdot}{\cdot})
\end{displaymath} is uniformly continuous.
\end{axiom}

In particular this means that if $m_\alpha(\phi_1,\phi_2)=0$ then
the two arrangements $\alpha\phi_1$ and $\alpha\phi_2$ are
identical.

\begin{axiom}
$m$ is upper semi-continuous with respect to the distribution
metric.  That is to say, for every $\varepsilon >0$, there exists
$\delta = \delta(\varepsilon,\alpha,\phi)$, such that if
$\distd{(\alpha,\phi)}{(\beta,\psi)} < \delta$ then $m(\beta,\psi)
< m(\alpha,\phi) + \varepsilon$.
\end{axiom}

This last condition implies that if the two measures
$(f^{\alpha,\phi})^*(\mu)$ and $(f^{\beta,\psi})^*(\nu)$ are the
same, then $m(\alpha,\phi)=m(\beta,\psi)$.  Hence the value $m$ is
well defined on those measures on $G^G$ which arise as such an
image, and we can write

\[m(\alpha,\phi)=m((f^{\alpha,\phi})^*(\mu)).\]

We can now define $m$-equivalence of two arrangements.
\begin{defn}  We say $\alpha$ and $\beta$
are $m$-equivalent arrangements if there exist $\phi_i$ which are
$m_\alpha$-Cauchy, $\phi_i^{-1}$ are $m_\beta$-Cauchy and
$\alpha\phi_i$ converges in probability to $\beta$.
\end{defn}

One can now define $m$-equivalence of actions on distinct spaces
as meaning there are conjugate versions of the actions as
arrangments on the same orbit space where the arrangements are
$m$-equivalent in the sense of the definition.

 We now give
a list of  examples to indicate the range of equivalence relations
that can be brought under this perspective.

Many examples of sizes have the common feature of being integrals
of some pointwise calculation of the distortion of a single orbit.
To make this precise we first review some material about
bijections of $G$.  Remember that $\mathcal B$ is the space of all
bijections of the group $G$ with the product topology,
$\mathcal{G}$ is the space of bijections fixing $\mathrm{id}$ and
we metrized both with a complete metric $d$.  The group $G$ can be
regarded as a subgroup of $\mathcal B$ acting by left
multiplication, ($g(g')=gg'$).  The map $\hat
H:\mathcal{B}\to\mathcal{G}$ given by $\hat
H(q)=qq(\mathrm{id})^{-1}$ is a contraction in $d$.  Also $G$
acting by right multiplication conjugates $\mathcal B$ to itself
giving an action of $G$ on $\mathcal B$.
($T_g(q)(g')=q(g'g)g^{-1}.$) We view this action by representing
an element $q\in\mathcal B$ by a map $f:G\to G$,
$f(g)=q(g)g^{-1}$.  Those maps $f\in G^G$ that arise from
bijections are a $G_\delta$ and hence a Polish space we call $F$.
The map $q\to f$ is obviously a homeomorphism from $\mathcal B$ to
$F$. For $f\in F$ let $Q(f)$ be the associated bijection and for
$q\in\mathcal{B}$ let $F(q)$ be the associated name in $G^G$. The
action of $G$ on $\mathcal B$ in its representation as $F$ is the
shift action $\sigma_g(f)(g')=f(g'g).$  Any rearrangement pair
$(\alpha,\phi)$ then gives rise to an ergodic shift invariant
measure on this Polish subset of $G^G$ and any ergodic shift
invariant measure is an ergodic action of $G$ with a canonical
rearrangement pair.   The probability measures on a Polish space
are weak* Polish and hence the invariant and ergodic measures on
this Polish space are weak* Polish.

We will now define a general class of sizes that arise as
integrals of valuations made on the bijections
$q^{\alpha,\phi}_x$.

\begin{defn}\index{size kernel}
  A Borel $D:\mathcal{B}\to{\mathbb R}^+$ is called a {\bf size kernel}
if it satisfies:
\begin{enumerate}
\item $D(q)\geq 0.$
\item $D(\mathrm{id})=0.$
\item $D(q(\mathrm{id})^{-1}q^{-1}q(\mathrm{id}))=D(q).$
\item $D(q_1(\mathrm{id})q_2q_1^{-1}(\mathrm{id}) q_1)\leq D(q_1)+D(q_2).$
\item For every $\varepsilon>0$ there is a
$\delta>0$ so that if $D(q)<\delta$ then
$d({\mathrm{id}},H(q))<\varepsilon.$
\item The function $\mu\to \int D(q(f))\,d\mu$ is weak*  continuous on the
space of shift invariant measures $\mu$ on the Polish space $F$.
\end{enumerate}
{\rm {\bf Note.} An element of $G$ is regarded as an element of
$\mathcal B$ acts by left multiplication.}
\end{defn}

For a size kernel $D$ we define
\[m^D(\alpha,\phi)=\int D(q_x^{\alpha,\phi})\,d\mu(x).\]
We call such an $m^D$ an {\bf integral size.}

\begin{example}[Conjugacy and Orbit Equivalence]\end{example}

These first two examples are the extremes of what is possible.
For one the equivalence class will be simply the full group orbit
and for the other it will be the entire set of arrangements.
  Both of the pseudometrics $d(q,\mathrm{id})$ and $d(H(q),\mathrm{id})$ are easily seen
to be size kernels and so both
\begin{eqnarray*}m^1(\alpha,\phi)&=\|(\alpha,\phi),(\alpha,\mathrm{id})\|_\alpha^s\text{ and}\\
  m^0(\alpha,\phi)&=\|(\alpha,\phi),(\alpha,\mathrm{id})\|_\alpha^w
\end{eqnarray*}
are sizes.

As $d$ makes $\mathcal B$ complete, relative to $m^1$ sequence
$phi_i$ is $m_\alpha^1$ Cauchy iff $\phi_i\to\phi$ in probability.
Thus $\alpha\overset{m^1}{\sim}\beta$ iff $\beta=\alpha\phi$ i.e.
they differ by an element of the full group and the equivalence
class of $\alpha$ is exactly its full group orbit. As for $m^0$,
for any $\alpha$ and $\beta$ one can  construct a sequence of
$\phi_i$ with $\alpha\phi_i\to\beta$ in $L^1$ with the sequence
$\phi_i$ an $m_\alpha$ Cauchy sequence.  Thus all arrangements are
$m^0$ equivalent.

\begin{example}[Kakutani Equivalence]\end{example}

For this example let $G=\mathbb{Z}^n$ and $B_N=[-N,N]^n$ be the
standard F\o lner sequence of boxes centered at $\vec 0$.  We
begin with a metric on $\mathbb{Z}^n$ given by
\[\tau(\vec u,\vec v)=\min\bigl(\|(\vec u/\|\vec u\|)-(\vec v/\|\vec v\|)\| +\bigl|\ln(\|\vec v\|)-\ln(\|\vec u\|)\bigr|,1\bigr)\]
(assuming $\vec 0/\|\vec 0\|=\vec 0$). What is important about
$\tau$ are the following two properties:
\begin{enumerate}
  \item $\tau$ is a metric on $\mathbb{Z}^n$ bounded by 1 and
\item $\vec u$ and $\vec v$ are $\tau$ close iff the norm of their difference is small in proportion to both of their norms.
\end{enumerate}

For $h\in\mathcal{G}$ set $B_N(h)=\{\vec v\in B_N | h(\vec v)\in
B_N\}$  (those elements of $B_N$ mapped into $B_N$ by $h$).  Now
set
\[k(h)=\sup_N\Bigl(\frac 1{\# B_N}(\sum_{\vec v\in B_N(h)}\tau(\vec v, h(\vec v))+\#\{\vec v\in B_N|h(\vec v)\notin B_N\})\Bigr).\]
Now set $K(q)=k(H(q).$
\begin{lem} The function $K$ is a size kernel.

\end{lem}

For $d=1$ standard arguments imply that this size yields even
Kakutani equivalence.  For $d>1$ it is leads to an analogous
equivalence relation among Katok cross-sections of $\mathbb
R^d$-actions.

Our last example moves beyond size kernels.

\begin{example}[Entropy as a Size]\end{example}

We discuss this example only for actions of $\mathbb Z$ although
the ideas extend to general countable amenable groups.

The size at its base will simply be the entropy of the
rearrangment itself.  We make this precise as follows. The
function $g_{(\alpha,\phi)}(x)=\alpha(x,\phi(x))$ takes on
countably many values and hence can be regarded as a countable
partition $g_{(\alpha,\phi)}$ of $X$.  Set $\Gamma_0^\alpha$ to be
those $\phi$ for which $g_{(\alpha,\phi)}$ is finite.  It is not
difficult to see that $\Gamma_0^\alpha$ is a subgroup and moreover
$\Gamma_0^{\alpha\psi}=\psi^{-1}\Gamma_0^\alpha\psi$ as
$g_{(\alpha\psi,\psi^{-1}\phi\psi)}(\psi^{-1}(x))=g_{(\alpha,\phi)}(x).$
It can be shown that the $\Gamma_0^\alpha$ are all $m^1_\alpha$
dense in $\Gamma$. For $\phi\in\Gamma_0^\alpha$ one can use the
entropy of the process $h(T^{\alpha},g_{(\alpha,\phi)})$ to start
the definition of a size defining
\[e(\alpha,\phi)=\inf_{\phi'\in\Gamma_0^\alpha} h(T^\alpha,g_{(\alpha,\phi')})+\mu\{x|\phi(x)\neq \phi'(x)\}.\]
Now set the size to be
\[m^e(\alpha,\phi)=e(\alpha,\phi)+m^0(\alpha,\phi).\]

\begin{prop}
  Two $\mathbb Z$-actions are $m^e$ equivalent iff they have the same entropy.
\end{prop}

\section*{2. Transference via orbit equivalence} \addsec
\vskip-5mm \hspace{5mm}

A second natural type of restriction can be placed on an orbit
equivalence.  Here the interest is in two arrangements $\alpha$
and $\beta$ of perhaps distinct groups.  Suppose $\mathcal A$ is
an invariant sub $\sigma$-algebra for the action $T^\alpha$ of
$G_1$ and $\beta$ is a $G_2$ arrangement of the same orbit space.
We say the orbit equivalence from $\alpha$ to $\beta$ is $\mathcal
A$-measurable if the function $h^{\alpha,\beta}(x,g_1)=\beta(x,
T^\alpha_{g_1}(x)$ describing or the orbit is rearranged, is
$\mathcal A$ measurable for all choices of $g_1$.  Up to conjugacy
we can regard all ergodic actions of infinite discrete and
amenable groups as residing on the same orbit space, so beyond
ergodicity no dynamical property will be preserved by orbit
equivalence. On the other hand, many dynamical properties have
versions ``relative to" an invariant sub $\sigma$-algebra and many
such properties are indeed invariant under orbit equivalences that
are measurable with respect to that sub $\sigma$-algebra.

This method first arose in work with B. Weiss showing that actions
of discrete amenable groups that have completely positive entropy
(cpe), commonly called $K$-systems, are mixing of all orders.
This transference method has been applied to a variety of
questions.  Here is an outline of the argument to this first
result to exhibit the format.  The complete argument can be found
in {\bf Entropy and mixing for amenable group actions} by D.J.
Rudolph and B. Weiss, {\it Annals of Mathematics}, {\bf 151},
(2000)m pp. 1119-1150.
\begin{lem} If $T$ is an action of a discrete amenable group $G$ and $T\times B$, its direct product with a Bernoulli
  action $B$ of $G$, is relatively cpe with respect to the Bernoulli second coordinate, then $T$ must be cpe.
\end{lem}
\begin{lem} If $\hat T$ and $\hat S$ are ergodic actions on the same orbits of two discrete amenable
  groups $G_1$ and $G_2$ and the orbit equivalence between them is $\mathcal A$ measurable where
$\mathcal A$ is a $\hat T$ invariant sub $\sigma$-algebra, then
for any partition $P$, the conditional entropies $h(T,P|{\mathcal
A})$ and $h(S,P|{\mathcal A})$ are equal.
\end{lem}

Let $S_i$ be a list of finite subsets of either $G_i$.  We say
the $S_i$ {\bf spread} if any particular $\gamma\neq {\mathrm
id}$ belongs to at most finitely many of the sets $S_iS_i^{-1}$.
If the sets $S_i(x)$ are random sequences of finite sets
depending on $x$, we can again say they are spread if for a.e.
$x$ they form a spread sequence.  A classical characterization of
the $K$-systems, which we state in a relative form says
\begin{thm}  $T$, a $\mathbb Z$-action, is relatively cpe with respect to a sub $\sigma$-algebra $\mathcal A$
  iff for all partitions $P$ and all $\mathcal A$ measurable and spread random
sequences of sets $S_i(x)$,
         \[\frac 1{\# S_i}\bigl[h(\underset {g\in S_i}\vee T_g(P)|{\mathcal A}) -\sum_{g\in S_i}h(T_g(P)|{\mathcal A})\bigr]
\underset i \to 0\] in $L^1$.
\end{thm}

That is to say, the translates of $P$ become conditionally ever
more independent the more spread they become.  Refer to this
property as $\mathcal A$-relative uniform mixing if it holds for
all $P$.
\begin{lem}
If $T$ and $S$ are $\mathcal A$-measurably orbit equivalent
actions of perhaps distinct groups and $T$ is $\mathcal
A$-relatively uniformly mixing.  Then $S$ is also.
 \end{lem}

We now describe the orbit transference proof that cpe actions of
discrete amenable groups are always mixing of all orders.  Suppose
$T$ is a cpe action of the group $G$.  Take $T\times B$ where $B$
is a Bernoulli action of $G$.  This direct product will be ergodic
and in fact relatively cpe with respect to the Bernoulli
coordinate.  Now $B$ is orbit equivalent to a $\mathbb Z$ action
and this orbit equivalence lifts to an orbit equivalence of
$T\times B$ to some ergodic $\mathbb Z$ action $S$.  The orbit
equivalence will be $\mathcal A$ measurable where $\mathcal A$ is
this Bernoulli coordinate algebra.  Now $S$ will still be
$\mathcal A$ relatively cpe and hence $\mathcal A$ relatively
uniformly mixing.  But now this tells us $T\times B$ is also
$\mathcal A$ relatively uniformly mixing.  Restricting this to
partitions that are measurable with respect to the first
coordinate tells us $T$ itself is uniformly mixing (without any
conditioning) and hence mixing of all orders.

A second and quite significant application of this method, due to
Dooley and Golodets, is to show that cpe actions have countable
Haar spectrum.  In as yet unwritten work, again with B. Weiss, one
can show that weakly mixing isometric extensions of Bernoulli
actions must be Bernoulli.

This remains an area of very active work. We end on an open
question.  Consider the known result for $\mathbb Z$ actions, that
a weakly mixing and isometric extension of a base action that is
mixing must itself be mixing.  Is this result true for general
amenable group actions?  To apply the transference method one
needs a relativized version of the result for $\mathbb Z$ actions.
That is to say, one needs to know that a relatively weakly mixing
relatively isometric extension of a relatively mixing action is
still relatively mixing.  What seems an obstacle here is simply
the definition of relative mixing over a sub $\sigma$-algebra
$\mathcal A$.  Certainly it means that for any sets $A$ and $B$
that
\[
  \bigl |E(I_A  I_B\circ T^j|\mathcal A)- E(I_A|\mathcal A)E(I_B\circ T^j|\mathcal
A)\bigr |\to 0 . \]
The question is, in what sense should it tend
to zero.  Pointwise convergence behaves well for orbit equivalence
but the above relativized question seems answerable only for mean
convergence.

\label{lastpage}

\end{document}